\nopagenumbers
\input amstex \magnification=\magstep 1 \baselineskip=21pt

\centerline {\bf La Derivada del Coseno}

El c\'\i rculo 
es $x^2 + y^2 = 1$. Al tomar la
diferencial de ambos lados, obtenemos\vskip-.2in $$2xdx + 2ydy = 0.$$
\vskip -.06in \noindent Al resolver para $dy$  obtenemos $dy = \frac{-x}y dx$.  Por 
lo tanto, la pendiente de la recta tangente es igual a $\frac{-x}y$.  ?`Cual es el 
negativo reciproco de esto?  Es la pendiente de la recta normal, $\frac yx$ pero 
es igual a la pendiente de la linea recta del origen al punto $(x,y)$, a saber, 
el radio.  Entonces, 
la recta tangente es normal al radio (Euclides, Lib.\ III Prop.\ 18).

\vskip -.06in	Ahora bien, $y = \text{sen\ }
\hskip-1.7pt\theta$ por definici\'on, y $x = \cos\theta$   por definici\'on, 
donde 
 $\theta$  es el largo del arco desde el tres en la esfera de reloj a $(x,y)$. Nos gustar\'\i a hallar $dx\over d\theta$. 

	F\'\i jese que $d\theta$ es el cambio en la largura del arco por la recta tangente al c\'\i rculo en $(x,y)$.  As\'\i\ que hay un tri\'angulo diferencial con una base $dx$, altura $dy$, e hipotenusa $d\theta$.
Los dos \'angulos agudos los llamemos $\alpha$ y $\frac\pi 2 - \alpha$.

\vskip 1.5in

  \includegraphics{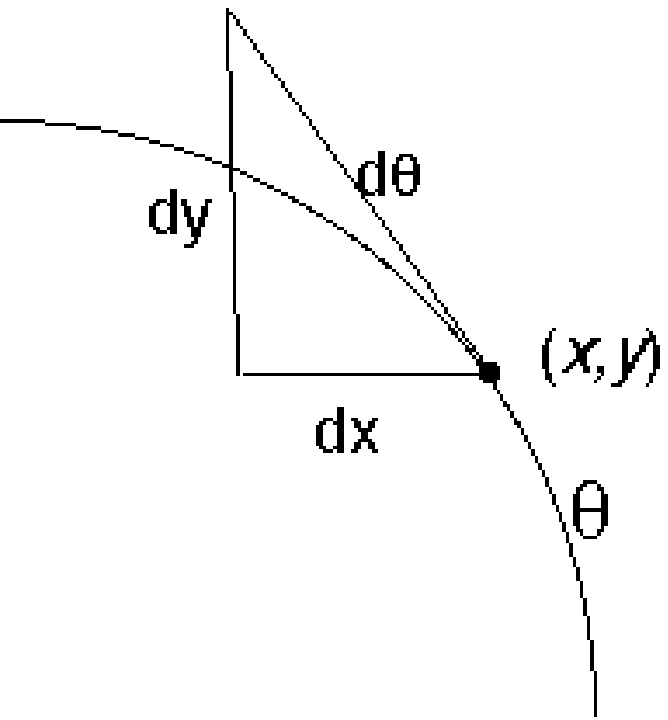} 

\vskip -.55in\quad\quad\quad\quad\quad\quad\quad\quad\quad\quad\quad\quad\quad\quad\quad\quad\quad\quad\quad\quad\quad\includegraphics{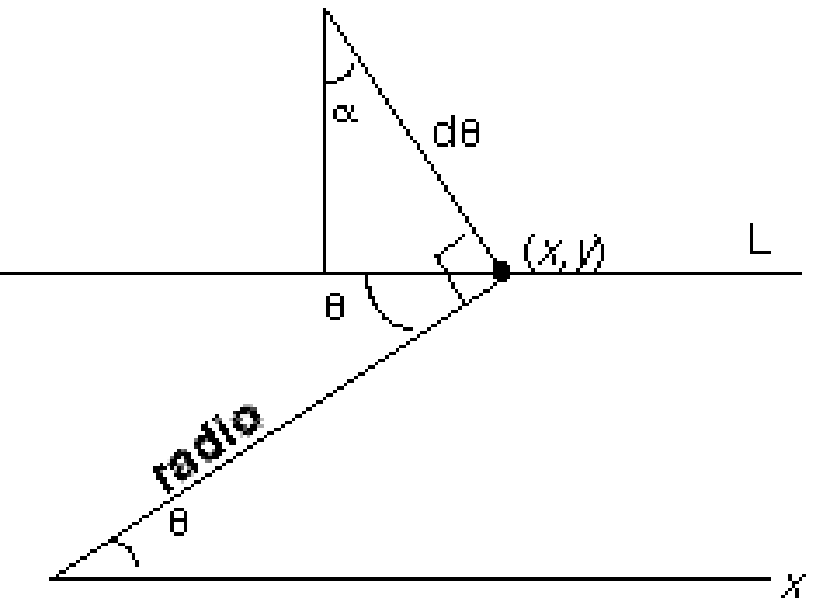} 

\noindent Pero la recta tangente ($d\theta$) es normal al radio.
  Entonces, el \'angulo marcado es recto.

\vskip -.06in El \'angulo \it alpha\rm\ del tri\'angulo, ?`es equivalente a 
$\theta$?  
Porque la linea $L$ (que extiende la base del tri\'angulo) es 
paralela al eje-$x$ , Euclides I, 29, dice que los dos \'angulos 
marcados como $\theta$ son bien iguales.
Pero el \'angulo $\theta$ es el complemento del \'angulo $\frac\pi 2 -\alpha$ del 
tri\'angulo.  
Por lo tanto, $\theta = \frac\pi 2 -(\frac\pi 2 -\alpha)=\alpha.$ 

\vskip -.06in Ahora bien el seno de $\alpha$\ es igual al cociente de la base
(que est\'a al lado opuesto) 
por la hipotenusa, \it esto es \rm, $ {-dx \over d\theta}$\ (la longitud de la base 
es $-dx$ porque $dx$ es negativo).  
$$\text{Entonces,}\hskip .6in{d\cos \theta \over d\theta}= -\hskip.4pt\text{sen\ }\hskip-2pt\theta,  \text{\quad\quad Q.E.D.}$$
\end

El lado opuesto al \'angulo $\alpha$ esta la base 
(que esta $-dx$) 
y por tanto el seno de 
$\theta$ es igual a la cocienta de la base para la hipotenusa....

The tangent is normal to the radius, by Euclid, ...
and so the marked angle is a right angle.  

The acute angles of the infinitesimal triangle are \alpha
(marked) and $\frac\pi 2 - \alpha$.

Because the line L (extending the base of the triangle) is 
parallel to the x-axis, Euclid ....
says that the angles marked $\theta$ are iguales.
Porque la linea $L$ (que extendar la base del tri\'angulo) es 
paralela al eje-$x$, Euclides dice que los dos \'angulos 
marcados como $\theta$ son bien iguales.

Now $\theta$ and $\frac\pi 2 - \alpha$ together comprise 
the right angle between the tangent and the radius, 
therefore they add up to $\frac\pi 2$ and so $\theta=\alpha$.

Therefore $\sin\theta = $ side opposite, which is $-dx$, 
divided by the hypotenuse, which is $d\theta.$